\documentclass[12pt]{article}
\usepackage{graphicx}
\usepackage{amsmath,amsthm,amssymb,enumerate}
\usepackage{euscript,mathrsfs}
\usepackage{color}
\usepackage{dsfont}
\usepackage{url}
\usepackage[left=2cm,right=2cm,top=3.5cm,bottom=3.5cm]{geometry}
\usepackage{color}
\usepackage[framemethod=tikz]{mdframed}
\allowdisplaybreaks

\usepackage{soul}

\catcode`\@=11 \@addtoreset{equation}{section}

\catcode`\@=12

\newtheorem{Theorem}{Theorem}[section]
\newtheorem{Proposition}[Theorem]{Proposition}
\newtheorem{Lemma}[Theorem]{Lemma}
\newtheorem{Corollary}[Theorem]{Corollary}

\theoremstyle{definition}
\newtheorem{Definition}[Theorem]{Definition}

\newtheorem{Remark}[Theorem]{Remark}

\newcommand{\bTheorem}[1]{
	\begin{Theorem} \label{T#1} }
	\newcommand{\eT}{\end{Theorem}}

\newcommand{\bProposition}[1]{
	\begin{Proposition} \label{P#1}}
	\newcommand{\eP}{\end{Proposition}}

\newcommand{\bLemma}[1]{
	\begin{Lemma} \label{L#1} }
	\newcommand{\eL}{\end{Lemma}}

\newcommand{\bCorollary}[1]{
	\begin{Corollary} \label{C#1} }
	\newcommand{\eC}{\end{Corollary}}

\newcommand{\bRemark}[1]{
	\begin{Remark} \label{R#1} }
	\newcommand{\eR}{\end{Remark}}

\newcommand{\bDefinition}[1]{
	\begin{Definition} \label{D#1} }
	\newcommand{\eD}{\end{Definition}}

\newcommand{\Del}{\Delta_x}

\newcommand{\Ds}{\mathbb{D}_x}

\newcommand{\vme}{\vm_\ep}

\newcommand{\bfomega}{\boldsymbol{\omega}}

\newcommand{\tvm}{\widetilde{\vc{m}}}
\newcommand{\tS}{\widetilde{S}}
\newcommand{\bfphi}{\boldsymbol{\varphi}}

\newcommand{\bFormula}[1]{
	\begin{equation} \label{#1}}
	\newcommand{\eF}{\end{equation}}

\newcommand{\Ov}[1]{\overline{#1}}

\newcommand{\Curl}{{\bf curl}_x}
\newcommand{\DC}{C^\infty_c}

\newcommand{\aleq}{\stackrel{<}{\sim}}

\newcommand{\vr}{\varrho}
\newcommand{\vre}{\vr_\ep}

\newcommand{\vte}{\vt_\ep}
\newcommand{\vue}{\vu_\ep}
\newcommand{\tvr}{\wtilde \vr}
\newcommand{\tvu}{{\wtilde \vu}}
\newcommand{\tvt}{\wtilde \vt}

\newcommand{\vt}{\vartheta}
\newcommand{\vu}{\vc{u}}
\newcommand{\vm}{\vc{m}}

\newcommand{\vc}[1]{{\bf #1}}

\newcommand{\Div}{{\rm div}_x}
\newcommand{\Grad}{\nabla_x}

\newcommand{\dx}{\,{\rm d} {x}}

\newcommand{\dt}{\,{\rm d} t }

\newcommand{\vU}{\vc{U}}

\newcommand{\intO}[1]{\int_{\Omega} #1 \ \dx}

\newcommand{\intOe}[1]{\int_{\Omega_\ep} #1 \ \dx}

\newcommand{\D}{{\rm d}}

\newcommand{\ep}{\varepsilon}

\newcommand{\R}{\mathbb{R}}
\newcommand{\vtB}{\vt_B}

\newcommand{\br}{ \nonumber \\ }

\def\softd{{\leavevmode\setbox1=\hbox{d}%
		\hbox to 1.05\wd1{d\kern-0.4ex{\char039}\hss}}}
\definecolor{Cgrey}{rgb}{0.85,0.85,0.85}
\definecolor{Cblue}{rgb}{0.50,0.85,0.85}
\definecolor{Cred}{rgb}{1,0,0}
\definecolor{fancy}{rgb}{0.10,0.85,0.10}
\definecolor{amaranth}{rgb}{0.9, 0.17, 0.31}

\newcommand\Cbox[2]{%
	\newbox\contentbox%
	\newbox\bkgdbox%
	\setbox\contentbox\hbox to \hsize{%
		\vtop{
			\kern\columnsep
			\hbox to \hsize{%
				\kern\columnsep%
				\advance\hsize by -2\columnsep%
				\setlength{\textwidth}{\hsize}%
				\vbox{
					\parskip=\baselineskip
					\parindent=0bp
					#2
				}%
				\kern\columnsep%
			}%
			\kern\columnsep%
		}%
	}%
	\setbox\bkgdbox\vbox{
		\color{#1}
		\hrule width  \wd\contentbox %
		height \ht\contentbox %
		depth  \dp\contentbox
		\color{black}
	}%
	\wd\bkgdbox=0bp%
	\vbox{\hbox to \hsize{\box\bkgdbox\box\contentbox}}%
	\vskip\baselineskip%
}

\mdfdefinestyle{MyFrame}{%
	linecolor=black,
	outerlinewidth=1pt,
	roundcorner=5pt,
	innertopmargin=\baselineskip,
	innerbottommargin=\baselineskip,
	innerrightmargin=10pt,
	innerleftmargin=10pt,
	backgroundcolor=white!20!white}

\newcommand{\tsl}{\textsl}

\newcommand{\veps}{\varepsilon}

\newcommand{\wtilde}{\widetilde}

\newcommand{\g}{\gamma}

\newcommand{\de}{\delta}
\renewcommand{\o}{\omega}

\newcommand{\bt}{\beta}

\allowdisplaybreaks

\begin{document}


\title{\bf Thermally driven fluid convection in the incompressible limit regime}

\author{Francesco Fanelli 
\thanks{The work of F.F. has been partially supported by the project CRISIS (ANR-20-CE40-0020-01), operated by the French National Research Agency (ANR).}
\and
Eduard Feireisl \thanks{The work of E.F. was partially supported by the
		Czech Sciences Foundation (GA\v CR), Grant Agreement
		21--02411S. The Institute of Mathematics of the Academy of Sciences of
		the Czech Republic is supported by RVO:67985840. } }

\date{}

\maketitle

{\small

\centerline{BCAM -- Basque Center for Applied Mathematics}
\centerline{Alameda de Mazarredo 14, E-48009 Bilbao, Basque Country, Spain}
\medbreak
\centerline{Ikerbasque -- Basque Foundation for Science}
\centerline{Plaza Euskadi 5, E-48009 Bilbao, Basque Country, Spain}
\medbreak
\centerline{Univ. Lyon, Universit\'e Claude Bernard Lyon 1, CNRS UMR 5208, Institut Camille Jordan,}
\centerline{F-69622 Villeurbanne, France} 
\medbreak
\centerline{Email address: \texttt{ffanelli@bcamath.org}}

\bigbreak
\bigbreak
\centerline{Institute of Mathematics of the Academy of Sciences of the Czech Republic}
\centerline{\v Zitn\'a 25, CZ-115 67 Praha 1, Czech Republic}
\medbreak
\centerline{Email address: \texttt{feireisl@math.cas.cz}}

}

\medbreak

\begin{abstract}
	
We consider a scaled Navier--Stokes--Fourier system describing the motion of a compressible, heat-conducting, viscous fluid driven 
by inhomogeneous boundary temperature distribution together with the gravitational force of a massive object placed outside the fluid. 
We identify the limit system in the low Mach/low Froude number regime for the ill prepared initial data. 
The fluid is confined to a bounded cavity with acoustically hard boundary enhancing reflection of acoustic waves.

\end{abstract}


{\small

\noindent
{\bf 2020 Mathematics Subject Classification:}{ 35B40 
(primary);
35Q30, 76N15, 35B25
(secondary).}

\medbreak

\tableofcontents

}

\section{Introduction}
\label{i}

Thermally driven convective flows are ubiquitous in many real world phenomena, notably in geophysics, meteorology and stellar dynamics or, in general, plasma 
motion, to name only a few. The dynamics is enhanced by mutual interaction of thermal expansion of a 
fluid with other volume forces, typically the gravitation.

Apparently, the fluid mass density is changing during the process and a correct mathematical description is therefore based on models of \emph{compressible} 
fluids. Still the simplified \emph{incompressible} formulation is frequently used notably in numerical simulations, see Barletta et al. \cite{barletta2022Bous}, \cite{barletta2022use},  Fr\" ohlich \cite{FrLaPe}, or the survey by Klein et al. \cite{KBSMRMHS}.
A natural idea is to justify the incompressible \emph{target} problem as a singular limit of the complete \emph{primitive} system. In the context of compressible viscous fluids, the primitive model is represented by the Navier--Stokes--Fourier 
system specified below, while the \emph{target} problem is 
usually
identified with the Oberbeck--Boussinesq approximation.

The limit process is described in detail 
by Zeytounian \cite{ZEY}, and the rigorous justification obtained in \cite{FN5}, \cite[Chapter 5]{FeNo6A} in the case of \emph{conservative} boundary conditions.
A suitable theoretical platform based on a new concept of weak solution to attack problems involving thermally driven fluid convection has been 
developed only recently in \cite{ChauFei}, \cite{FeiNovOpen}. The singular limit problem in this new framework has been then revisited in \cite{BelFeiOsch}. Rather surprisingly and in contrast 
with the previously obtained results, thermally driven flows give rise to an additional ``unexpected'' non--local term in the boundary conditions for the target 
problem, see \cite[Theorem 4.1]{BelFeiOsch}. The result \cite{BelFeiOsch} relies on the relative energy method and
yields convergence to the target system in a strong topology associated to the energy norm. In particular, the convergence rate can be estimated
in terms of the initial data. However, this approach imposes two major restrictions:
\begin{itemize}
 \item \emph{well--prepared} initial data for the primitive system;
 \item the existence of smooth solutions for the target system, in particular the convergence result can be stated
 only on a short time interval in the natural 3--d setting, cf. \cite{Abb-Feir}.
\end{itemize}

Our goal is to consider the \emph{ill--prepared} initial data
to show convergence to a global--in--time solution of the target problem in the framework of weak topology.
On the one hand, the ill--prepared data allow  for relatively large deviations from the equilibrium state.
On the other hand, thermally driven systems operate in the far from equilibrium regime,
as observed \tsl{e.g.} in the well known Rayleigh--B\'enard convection model.
In addition, analogously to the conservative case (Neumann boundary conditions for the temperature), for which we refer to \cite[Section 5.5]{FN5},
the presence of the acoustically hard boundary enhances oscillations of acoustic waves.
In particular, the convergence 
to the target system can be shown only in the weak Lebesgue topology,
which entails non--trivial difficulties completely absent whenever the initial data are well--prepared. The benefits obtained by this approach are:

\begin{itemize}
	\item a large class of initial data allowing the system to stay out of thermodynamics equilibrium; 
	
	\item convergence to global--in--time weak solutions of the target problem.
	\end{itemize}

\noindent
More details will be discussed in Section \ref{ss:idea-proof} below.

\subsection{Primitive Navier--Stokes--Fourier system} \label{ss:primitive}

We adopt the \emph{Navier--Stokes--Fourier (NSF) system } as the exact -- primitive -- model. The state of the fluid is characterized 
by its mass density $\vr = \vr(t,x)$, the absolute temperature $\vt = \vt(t,x)$, and the velocity $\vu = \vu(t,x)$. Assuming
the Mach number  is proportional to a small parameter $\veps>0$, the
scaled system of field  equations written in the Eulerian coordinate frame reads:

\begin{mdframed}[style=MyFrame] 
\begin{align} 
	\partial_t \vr + \Div (\vr \vu) &= 0, \label{i1}\\
	\partial_t (\vr \vu) + \Div (\vr \vu \otimes \vu) + \frac{1}{\ep^2} \Grad p(\vr, \vt) &= \Div \mathbb{S}(\vt, \Ds \vu) + \frac{1}{\ep} \vr \Grad G, \label{i2} \\ 
	\partial_t (\vr s(\vr, \vt)) + \Div (\vr s (\vr, \vt) \vu) + \Div \left( \frac{ \vc{q} (\vt, \Grad \vt) }{\vt} \right) &= 
	\frac{1}{\vt} \left( \ep^2 \mathbb{S} : \Ds \vu - \frac{\vc{q} (\vt, \Grad \vt) \cdot \Grad \vt }{\vt} \right).
	\label{i3}	
\end{align}
\end{mdframed}
The quantity $s = s(\vr, \vt)$ in \eqref{i3} is the entropy of the system, related to the pressure $p = p(\vr, \vt)$ and the internal energy $e = e(\vr, \vt)$ through Gibbs' equation 
\begin{equation} \label{i10} 
	\vt D s = D e + p D \left( \frac{1}{\vr} \right).
\end{equation}
The viscous stress tensor is given by Newton's rheological law 
\begin{equation} \label{i7}
	\mathbb{S}(\vt, \Ds \vu) = 2\mu(\vt) \left( \Ds\vu - \frac{1}{3} \Div \vu \mathbb{I} \right) + \eta(\vt) \Div \vu \mathbb{I}, 
\end{equation}
where $\Ds\vu = \frac{1}{2}(\Grad \vu + \Grad^t \vu)$ is the symmetric part of the velocity gradient $\Grad \vu$. 
The internal energy flux is given by Fourier's law
\begin{equation} \label{i8} 
	\vc{q}(\vt, \Grad \vt) = - \kappa (\vt) \Grad \vt.
\end{equation}
The potential $G = G(x)$ represents the gravitational force imposed by a massive object placed outside the fluid cavity.

Besides the Mach number proportional to $\ep$, there is another scaled quantity, namely $\frac{1}{\ep} \vr \Grad G$, where $\sqrt{\ep}$ is the so--called Froude number. The fluid motion is therefore almost incompressible and mildly stratified, see \cite[Chapter 4]{FeNo6A}.

The fluid is confined to a \emph{bounded} physical domain $\Omega \subset R^3$ with acoustically hard boundary. Specifically, 
the velocity satisfies the complete slip boundary conditions,
\begin{align} 
	\vu \cdot \vc{n}|_{\partial \Omega} = 0,\qquad  
	\left( \mathbb{S}(\vt, \Ds \vu) \cdot \vc{n} \right) \times \vc{n}|_{\partial \Omega} = 0. 
	\label{i5}
	\end{align}
Finally, we impose appropriately scaled Dirichlet boundary conditions for the temperature, 	
\begin{equation} \label{i6}
	\vt|_{\partial \Omega} = \vtB \equiv \Ov{\vt} + \ep \mathfrak{T}_B, 
	\end{equation}
where $\Ov{\vt}$ is a positive constant. As we shall see below, condition \eqref{i6} corresponds to \emph{ill--prepared} boundary data for the temperature.

\subsection{Static solutions}
\label{S}

In accordance with the scaling of \eqref{i2}, \eqref{i3},
the zero--th order terms in the (formal) asymptotic limit
written in terms of powers of the small parameter $\ep$ is determined by the stationary (static) problem 
\begin{equation} \label{S1}
	\Grad p(\tvr_{\ep}, \Ov{\vt} + \ep \mathfrak{T}_B ) =  \ep \tvr_{\ep} \Grad G,
\end{equation}
where
\begin{equation} \label{T_B-harm}
\text{$\mathfrak{T}_B$ is the unique harmonic extension of the boundary datum
inside $\Omega$.}
     \end{equation}

Observe that equation \eqref{i1} conserves  the total mass, meaning that
\[
 \frac{\rm d}{\dt}\intO{\vre(t,\cdot)} = 0.
\]
Therefore, given some positive number $\Ov\vr>0$, we supplement equation \eqref{i1} with the condition that
\[
 \forall\ \veps>0,\qquad \frac{1}{|\Omega|} \intO{\vre} = \Ov{\vr}.
\]
Then, the solution $\tvr_{\ep}$ to equation \eqref{S1} is normalized so that 
\begin{equation}
\label{S3}
\intO{ (\Ov \vr - \tvr_\ep) } = 0,
\end{equation}
which implies in particular that
\begin{equation} \label{S2}
	\intO{ (\vre(t, \cdot) - \tvr_{\ep})} = 0 \qquad \mbox{ for any } 
	\ t \geq 0.
\end{equation}

Finally, as $G = G(x)$ is given and smooth,  arguing as  in \cite{DS-F-S-WK}, one gathers that
\begin{equation} \label{S4}
	\| \tvr_\ep - \Ov{\vr} \|_{L^\infty(\Omega)} \lesssim \ep .
\end{equation}
Here and always hereafter, the symbol $a \lesssim b$ means there is a constant $c > 0$, \emph{independent} of the scaling parameter $\ep$, such that $a \leq cb$.

\subsection{Target Oberbeck--Boussinesq system} \label{ss:target}

The asymptotic limit for well--prepared data has been identified in \cite[Theorem 4.1]{BelFeiOsch}. Our goal is to show that the ill--prepared data give rise to convergence to the same limit in the weak topology. Specifically, we show
\begin{align} 
	\frac{\vre - \Ov{\vr}}{\ep} \to \mathfrak{R}, \qquad
	\frac{\vte - \Ov{\vt}}{\ep} \to \mathfrak{T}, \qquad
\vue \to \vU,
\nonumber
\end{align}	
in a suitable weak topology, where the limit quantities
\[
\mathfrak{R}, \quad \vU, \quad \and \Theta := \mathfrak{T} - \frac{\lambda(\Ov{\vr}, \Ov{\vt})}{|\Omega|} \intO{ \mathfrak{T}}
\]
solve the \emph{Oberbeck--Boussinesq (OB) system}:
\begin{mdframed}[style=MyFrame]
	\begin{align} 
		\Div \vU &= 0, \br	
		\Ov{\vr} \Big( \partial_t \vU + \vU \cdot \Grad \vU \Big) + \Grad \Pi &= \Div \mathbb{S}(\Ov{\vt}, \Grad \vU) +  \mathfrak{R} \Grad G, \br
		\Ov{\vr} c_p(\Ov{\vr}, \Ov{\vt} ) \Big( \partial_t \Theta + \vU \cdot \Grad \Theta \Big)	- 
		\Ov{\vr} \ \Ov{\vt} \alpha(\Ov{\vr}, \Ov{\vt} ) \vU \cdot \Grad G
		&= \kappa(\Ov{\vt}) \Del \Theta, \br
		\frac{\partial p(\Ov{\vr}, \Ov{\vt} ) }{\partial \vr} \Grad \mathfrak{R} + 
		\frac{\partial p(\Ov{\vr}, \Ov{\vt} ) }{\partial \vt} \Grad \Theta	 &= \Ov{\vr} \Grad G	
		\label{ObBs1}
	\end{align}
	with the complete slip boundary conditions for the velocity
	\begin{equation} \label{i18bis}
		\vU \cdot \vc{n}|_{\partial \Omega} = 0,\ (\mathbb{S}(\Ov{\vt}, \Ds \vU) \cdot \vc{n}) \times \vc{n}|_{\partial \Omega} = 0
	\end{equation}
	and a \emph{non--local} boundary condition for the temperature deviation
	\begin{equation} \label{i14B}
		\Theta|_{\partial \Omega} = \mathfrak{T}_B - \frac{\lambda(\Ov{\vr}, \Ov{\vt})}{1 - \lambda(\Ov{\vr}, \Ov{\vt})} \frac{1}{|\Omega|}\intO{ \Theta }.	
	\end{equation}
\end{mdframed}

\noindent Here, the material parameters are:
\begin{itemize}
 \item the thermal expansion coefficient
\[
\alpha(\Ov{\vr}, \Ov{\vt} ) \equiv \frac{1}{\Ov{\vr}}  \frac{\partial p(\Ov{\vr}, \Ov{\vt} ) }{\partial \vt} \left( \frac{\partial p(\Ov{\vr}, \Ov{\vt} ) }{\partial \vr} \right)^{-1};
\]
\item the specific heat and constant pressure
\[
 c_p (\Ov{\vr}, \Ov{\vt} ) \equiv \frac{\partial e(\Ov{\vr}, \Ov{\vt} ) }{\partial \vt}	+ \Ov{\vr}^{-1} \Ov{\vt}
	\alpha(\Ov{\vr}, \Ov{\vt} ) \frac{\partial p(\Ov{\vr}, \Ov{\vt} ) }{\partial \vt};
\]
\item the coefficient
\[
\lambda(\Ov{\vr}, \Ov{\vt}) = \frac{\Ov{\vt} \alpha (\Ov{\vr}, \Ov{\vt} ) }
{\Ov{\vr} c_p(\Ov{\vr}, \Ov{\vt}) } \frac{\partial p(\Ov{\vr}, \Ov{\vt} )}{\partial \vt} \ \in (0,1).
\] 
\end{itemize}
Our working hypothesis specified below asserts $\alpha > 0$, meaning, in particular, $\partial_\vt p > 0$. 

	

\subsection{Singular limit, main ideas of the proof} \label{ss:idea-proof}

Our goal is to identify the asymptotic limit of the family $(\vre, \vte, \vue)_{\ep > 0}$ of weak solutions to the scaled (NSF) system 
emanating from the initial data 
\begin{align} 
	\vre(0,x) &= \tvr_\ep + \ep \vr_{0,\ep}, \br
	\vte(0,x) & = \Ov{\vt} + \ep \vt_{0,\ep}, \br
	\vue(0,x) & =\vu_{0,\ep},
\label{wpd}	  
\end{align} 
where we assume that
\begin{equation} \label{wpd1}
\| \vr_{0, \ep} \|_{L^\infty(\Omega)} +  
\| \vt_{0, \ep} \|_{L^\infty(\Omega)} + \| \vu_{0, \ep} \|_{L^\infty(\Omega; R^3)} \lesssim 1.
\end{equation}
The bounds \eqref{wpd1} characterize the ill--prepared nature of the data and can be slightly relaxed to suitable $L^p-$norms. 

As announced in the preceding section, our goal is to show that the scaled quantities 
\[
\frac{\vre - \Ov{\vr}}{\ep},\ \ \frac{\vte - \Ov{\vt}}{\ep} \ \mbox{ and }\ \vue 
\]
converge weakly to a (weak) solution of the Oberbeck--Boussinesq system \eqref{ObBs1}, \eqref{i18bis}, \eqref{i14B}.
The precise result is stated in Section \ref{m}, see Theorem \ref{TM1}.
Although the target dynamics has already been identified in \cite{BelFeiOsch},
the strategy of the proof of the convergence in the case of ill--prepared data is rather different.

More precisely, the main steps to the proof of Theorem \ref{TM1} can be summarised as follows.

\begin{enumerate}[(i)]
\item	In Section \ref{L}, we introduce the class of weak solutions to the primitive (NSF) system relevant for the Dirichlet boundary conditions for the temperature.
In comparison with \cite{ChauFei}, \cite{FeiNovOpen}, the entropy production rate 
is explicitly specified and controlled by the available 
bounds uniform with respect to the scaling parameter $\ep > 0$.

\item In Section \ref{m}, we state our main result.

\item In Section \ref{e}, we derive the necessary uniform bounds independent of the scaling parameter $\ep>0$,
the main tool being the relative energy inequality.

\item  The weak convergence towards a weak solution of the target Oberbeck--Boussinesq system is performed in Section \ref{w}, 
except the convective term $\vue\otimes\vue$ appearing in \eqref{i2}.

\item Finally, in Section \ref{A} we tackle the proof of the convergence of the convective term to the expected limit $\vu\otimes\vu$.
In general, non--linear compositions do not commute with weak convergence; however, the special structure of the convective term
allows us to conclude.
	\end{enumerate}

Let us give more details about point (v) above.
In our case, the Dirichlet boundary conditions imposed of the temperature destroy the structure of the acoustic 
equation  obtained in the case of conservative boundary conditions in \cite{FN5}. The main reason is that 
the average of the temperature perturbation
\[
\intO{ \frac{\vte - \Ov{\vt}}{\ep} } 
\]
is no longer a conserved quantity.
This fact is, of course, completely harmless, in the case of well--prepared data, where acoustic waves are not created.
On the contrary, when considering ill--prepared initial data, such a problem cannot be circumvented. 
Therefore, for studying
propagation of acoustic waves, a cut--off must be performed in the neighbourhood of the boundary $\partial \Omega$. 
Now, this cut--off operation produces extra driving terms in the associated wave equation, which are ``out of scaling'' with respect to the low Mach number regime of equations
\eqref{i1}-\eqref{i3}. The situation is similar to the ones appearing in the multiscale problems considered in \cite{DS-F-S-WK, Fan}, and can be handled
%
by resorting to a compensation compactness argument, in the spirit of the so--called local method 
proposed by Lions and Masmoudi \cite{LIMA6, LIMA-JMPA}, cf. also the survey by Masmoudi \cite{MAS1}
(see \tsl{e.g.} \cite{Fan, Fa-Ga, Fa-Za, F-GV-G-N} and references therein for recent applications of that method in different contexts).

\section{Weak solutions to the primitive NSF system}
\label{L}

Our analysis is carried out in the framework of the class of weak solutions to the NSF system
introduced in \cite{ChauFei}, cf.~also \cite{FeiNovOpen}.

\begin{Definition}[{\bf Weak solution to the NSF system}] \label{DL1}
	We say that a trio $(\vr, \vt, \vu)$ is a \emph{weak solution} of the (scaled) NSF system \eqref{i1}--\eqref{i8}, with the boundary conditions \eqref{i5}, \eqref{i6},
	and the initial data
	\[
	\vr(0, \cdot) = \vr_0,\ \vr \vu (0, \cdot) = \vr_0 \vu_0,\ 
	\vr s(0, \cdot) = \vr_0 s(\vr_0, \vt_0),
	\]
if the following holds:	
\begin{itemize}
	
	\item the solution belongs to the {\bf regularity class}: 
	\begin{align}
		\vr &\in L^\infty(0,T; L^\gamma(\Omega)) \ \mbox{for some}\ \gamma > 1,\ \vr \geq 0 
		\ \mbox{a.a.~in}\ (0,T) \times \Omega, \br
		\vu &\in L^2(0,T; W^{1,2} (\Omega; R^3)), \ \vu \cdot \vc{n}|_{\partial \Omega} = 0, \br 
		\vt^{\beta/2} ,\ \log(\vt) &\in L^2(0,T; W^{1,2}(\Omega)) \ \mbox{for some}\ \beta \geq 2,\ 
		\vt > 0 \ \mbox{a.a.~in}\ (0,T) \times \Omega, \br
		(\vt - \vtB) &\in L^2(0,T; W^{1,2}_0 (\Omega));
		\label{Lw6}
	\end{align}
	 
	\item the {\bf equation of continuity} \eqref{i1} is satisfied in the sense of distributions, namely
	\begin{align} 
		\int_0^T \intO{ \Big[ \vr \partial_t \varphi + \vr \vu \cdot \Grad \varphi \Big] } \dt &=  - 
		\intO{ \vr_0 \varphi(0, \cdot) }
		\label{Lw4}
	\end{align}
	for any $\varphi \in C^1_c([0,T) \times \Ov{\Omega} )$;
	\item the {\bf momentum equation} \eqref{i2} is satisfied in the sense of distributions, specifically 
	\begin{align}
		\int_0^T &\intO{ \left[ \vr \vu \cdot \partial_t \bfphi + \vr \vu \otimes \vu : \Grad \bfphi +   
			\frac{1}{\ep^2} p(\vr, \vt) \Div \bfphi \right] } \dt \br &= \int_0^T \intO{ \left[ \mathbb{S}(\vt, \Ds \vu) : \Grad \bfphi - \frac{1}{\ep} \vr \Grad G \cdot \bfphi \right] } \dt - 
		\intO{ \vr_0 \vu_0 \cdot \bfphi (0, \cdot) }
		\label{Lw5}
	\end{align}	
	for any $\bfphi \in C^1_c([0, T) \times \Ov{\Omega}; R^3)$ such that 
	$\bfphi \cdot \vc{n}|_{\partial \Omega} = 0$;
	
	\item the {\bf entropy balance} \eqref{i3} is replaced by 
	\begin{align}
		- \int_0^T &\intO{ \left[ \vr s(\vr, \vt) \partial_t \varphi + \vr s (\vr ,\vt) \vu \cdot \Grad \varphi + \frac{\vc{q} (\vt, \Grad \vt )}{\vt} \cdot 
			\Grad \varphi \right] } \dt \br &= \int_0^T \int_{\Omega}{ \varphi \ \D \sigma(t,x)}  + \intO{ \vr_0 s(\vr_0, \vt_0) 
			\varphi (0, \cdot) } 
		\label{Lw7} 
	\end{align}
	for any $\varphi \in C^1_c([0, T) \times \Omega)$, where 
the	 entropy production rate $\sigma$ is a non--negative Radon measure, 
$\sigma \in \mathcal{M}^+([0,T] \times \Ov{\Omega})$, 
satisfying 
\begin{equation} \label{Lw7a}
\sigma \geq \frac{1}{\vt} \left( \ep^2 \mathbb{S}(\vt, \Ds \vu) : \Ds \vu - \frac{\vc{q}(\vt ,\Grad \vt): \Grad \vt}{\vt}\right);	
	\end{equation}
	\item  the {\bf ballistic energy balance}
	\begin{align}  
		- &\int_0^T \partial_t \psi	\intO{ \left[ \ep^2 \frac{1}{2} \vr |\vu|^2 + \vr e(\vr, \vt) - \vtB \vr s(\vr, \vt) \right] } \dt  \br &+ \int_0^T \int_{\Ov{\Omega}} \psi \vtB\  \D \sigma(t,x) 
		 \br
		&\leq 
		\int_0^T \psi \intO{ \left[ \ep \vr \vu \cdot \Grad G - 
			\vr s(\vr, \vt) \partial_t \vt_B - \vr s(\vr, \vt) \vu \cdot \Grad \vtB - \frac{\vc{q}(\vt, \Grad \vt)}{\vt} \cdot \Grad \vtB \right] } \dt \br 
		&+ \psi(0) \intO{ \left[ \frac{1}{2} \ep^2 \vr_0 |\vu_0|^2 + \vr_0 e(\vr_0, \vt_0) - \vtB(0, \cdot) \vr_0 s(\vr_0, \vt_0) \right] }
		\label{Lw8}
	\end{align}
	holds true for any $\psi \in C^1_c ([0, T))$, with $\psi \geq 0$, and any smooth extension of the boundary datum $\vtB$.
\end{itemize}
 
	\end{Definition}

Definition \ref{DL1}, in the spirit of \cite[Chapter 3]{FeNo6A}, provides more information than its	counterpart introduced in \cite{ChauFei}, \cite[Chapter 12]{FeiNovOpen}. Indeed the entropy dissipation rate 
$\sigma$ appears in both \eqref{Lw7}, \eqref{Lw8}; in particular, the ballistic energy dissipation in \eqref{Lw8} controls the entropy dissipation modulo a multiplicative 
factor related to the boundary temperature $\vtB$.	
It can be shown that the measure $\sigma$ is identified in the course of the existence proof, in particular all results obtained 
in \cite{FeiNovOpen} remain valid in the new framework. 

\subsection{Relative energy inequality}

In addition to Gibbs' equation \eqref{i10}, we impose the hypothesis of thermodynamic stability written in the form 
\begin{equation} \label{HTS}
	\frac{\partial p(\vr, \vt) }{\partial \vr } > 0,\quad 
	\frac{\partial e(\vr, \vt) }{\partial \vt } > 0 \qquad \mbox{ for all }\quad \vr, \vt > 0.
\end{equation}
Following \cite{ChauFei}, we introduce the scaled \emph{relative energy} 
\begin{align}
	E_\ep &\left( \vr, \vt, \vu \Big| \tvr , \tvt, \tvu \right) \br &= \frac{1}{2}\vr |\vu - \tvu|^2 + 
	\frac{1}{\ep^2} \left[ \vr e - \tvt \Big(\vr s - \tvr s(\tvr, \tvt) \Big)- 
	\Big( e(\tvr, \tvt) - \tvt s(\tvr, \tvt) + \frac{p(\tvr, \tvt)}{\tvr} \Big)
	(\vr - \tvr) - \tvr e (\tvr, \tvt) \right] .
	\nonumber
\end{align}
Now, the hypothesis of thermodynamic stability \eqref{HTS} can be equivalently rephrased as (strict) convexity 
of the total energy expressed with respect to the conservative entropy variables
\[
E_\ep \Big( \vr, S = \vr s(\vr, \vt), \vm = \vr \vu \Big) \equiv \frac{1}{2} \frac{|\vm|^2}{\vr} + 
\frac{1}{\ep^2} \vr e(\vr, S),
\]
whereas the relative energy can be written as
\begin{align}
E_\ep &\left( \vr, S, \vm \Big| \tvr , \tS, \tvm \right) = E_\ep(\vr, S, \vm) - \left< \partial_{\vr, S, \vm} E_\ep(\tvr, \tS, \tvm) ; (\vr - \tvr, S - \tS, \vm - \tvm) \right> - E_\ep(\tvr, \tS, \tvm).
\nonumber
\end{align}

Finally, as observed in \cite{ChauFei}, any weak solution of the NSF system in the sense of Definition \ref{DL1} satisfies the \emph{relative energy inequality} in the form: 
 \begin{align}
	&\left[ \intO{ E_\ep \left(\vr, \vt, \vu \Big| \tvr, \tvt, \tvu \right) } \right]_{t = 0}^{t = \tau}  
	+ \frac{1}{\ep^2}\int_0^\tau \int_{\Ov{\Omega}}
	\tvt \ \D \sigma(t,x) \br 
	&\leq - \frac{1}{\ep^2} \int_0^\tau \intO{ \left( \vr (s(\vr, \vt) - s(\tvr, \tvt)) \partial_t \tvt + \vr (s(\vr,\vt) - s(\tvr, \tvt)) \vu \cdot \Grad \tvt -
		\frac{\kappa (\vt) \Grad \vt}{\vt}  \cdot \Grad \tvt \right) } \dt \br 
	&- \int_0^\tau \intO{ \Big[ \vr (\vu - \tvu) \otimes (\vu - \tvu) + \frac{1}{\ep^2} p(\vr, \vt) \mathbb{I} - \mathbb{S}(\vt, \Ds \vu) \Big] : \Ds \tvu } \dt   \br
	&+ \int_0^\tau \intO{ \vr \left[ \frac{1}{\ep} \Grad G  - \partial_t \tvu - (\tvu \cdot \Grad) \tvu  \right] \cdot (\vu - \tvu) } \dt  \br 
	&+ \frac{1}{\ep^2} \int_0^\tau \intO{ \left[ \left( 1 - \frac{\vr}{\tvr} \right) \partial_t p(\tvr, \tvt) - \frac{\vr}{\tvr} \vu \cdot \Grad p(\tvr, \tvt) \right] } \dt
	\label{L4}
\end{align}
for a.a. $\tau > 0$ and any trio of continuously differentiable functions $(\tvr, \tvt, \tvu)$ satisfying
\begin{equation} \label{L5}
	\tvr > 0,\quad \tvt > 0,\quad \tvt|_{\partial \Omega} = \vtB, \quad \tvu \cdot \vc{n}|_{\partial \Omega} = 0.
\end{equation}

\subsection{Constitutive relations}
\label{CR}

The existence theory developed in \cite{ChauFei} requires certain restrictions to be imposed on the 
constitutive relations (state equations) similar to those introduced in the monograph
\cite[Chapters 1,2]{FeNo6A}. Specifically, the equation of state reads
\[
p(\vr, \vt) = p_{\rm m} (\vr, \vt) + p_{\rm rad}(\vt), 
\]
where $p_{\rm m}$ is the pressure of a general \emph{monoatomic} gas, 
\begin{equation} \label{con1}
	p_{\rm m} (\vr, \vt) = \frac{2}{3} \vr e_{\rm m}(\vr, \vt),
\end{equation}
enhanced by the radiation pressure 
\[
p_{\rm rad}(\vt) = \frac{a}{3} \vt^4,\qquad a > 0.
\]
As observed in \cite{DF1}, the radiation pressure prevents uncontrolled temperature oscillations 
in the (hypothetical) vacuum zones and, as such, is indispensable for the global existence theory to be valid. 
The pressure $p_m$ can be more general in the sense specified in \cite[Chapter 1, Section 1.4]{FeNo6A}.

Accordingly, the internal energy reads 
\[
e(\vr, \vt) = e_{\rm m}(\vr, \vt) + e_{\rm rad}(\vr, \vt),\qquad e_{\rm rad}(\vr, \vt) = \frac{a}{\vr} \vt^4.
\]

The specific form of \eqref{con1} is another issue to be discussed, obviously. In the context of gases, 
the natural candidate is provided by the standard Boyle--Mariotte law $p = \vr \vt$. Unfortunately, this \textsl{ansatz} 
fails to provide even the expected energy estimates as long as the boundary temperature is prescribed. To circumvent 
this difficulty, more physics must be taken into account, cf. \cite[Chapter 1]{FeNo6A}, as specified here below.

\begin{itemize}
	
	\item {\bf Gibbs' relation} together with \eqref{con1} yield 
	\[
	p_{\rm m} (\vr, \vt) = \vt^{\frac{5}{2}} P \left( \frac{\vr}{\vt^{\frac{3}{2}}  } \right)
	\]
	for a certain $P \in C^1[0,\infty)$.
	Consequently, 
	\begin{equation} \label{w9}
		p(\vr, \vt) = \vt^{\frac{5}{2}} P \left( \frac{\vr}{\vt^{\frac{3}{2}}  } \right) + \frac{a}{3} \vt^4,\quad
		e(\vr, \vt) = \frac{3}{2} \frac{\vt^{\frac{5}{2}} }{\vr} P \left( \frac{\vr}{\vt^{\frac{3}{2}}  } \right) + \frac{a}{\vr} \vt^4, \qquad a > 0.
	\end{equation}
	
	\item {\bf Hypothesis of thermodynamic stability} \eqref{HTS} expressed in terms of 
	$P$ gives rise to
	\begin{equation} \label{w10}
		P(0) = 0,\qquad P'(Z) > 0 \ \mbox{ for }\ Z \geq 0,\qquad \frac{ \frac{5}{3} P(Z) - P'(Z) Z }{Z} > 0 \ \mbox{ for }\ Z > 0.
	\end{equation} 	
	In particular, the function $Z \mapsto P(Z)/ Z^{\frac{5}{3}}$ is decreasing, and we suppose 
	\begin{equation} \label{w11}
		\lim_{Z \to \infty} \frac{ P(Z) }{Z^{\frac{5}{3}}} = p_\infty > 0.
	\end{equation}
	
	\item 
	The associated {\bf entropy} takes the form 
	\begin{equation} \label{w12}
		s(\vr, \vt) = s_{\rm m}(\vr, \vt) + s_{\rm rad}(\vr, \vt),\qquad s_{\rm m} (\vr, \vt) = \mathcal{S} \left( \frac{\vr}{\vt^{\frac{3}{2}} } \right),\qquad
		s_{\rm rad}(\vr, \vt) = \frac{4a}{3} \frac{\vt^3}{\vr}, 
	\end{equation}
	where one has
	\begin{equation} \label{w13}
		\mathcal{S}'(Z) = -\frac{3}{2} \frac{ \frac{5}{3} P(Z) - P'(Z) Z }{Z^2} < 0.
	\end{equation}
	Finally, we impose the {\bf Third law of thermodynamics}, cf.~Belgiorno \cite{BEL1}, \cite{BEL2}, requiring the entropy to vanish 
	when the absolute temperature approaches zero, namely
	\begin{equation} \label{w14}
		\lim_{Z \to \infty} \mathcal{S}(Z) = 0.
	\end{equation}
	
\end{itemize}

As for the transport coefficients, we suppose they are continuously differentiable functions of the 
temperature satisfying
\begin{align} 
	0 < \underline{\mu}(1 + \vt) &\leq \mu(\vt),\qquad |\mu'(\vt)| \leq \Ov{\mu}, \br 
	0 &\leq \underline{\eta} (1 + \vt) \leq \eta (\vt) \leq \Ov{\eta}(1 + \vt), \br
	0 < \underline{\kappa} (1 + \vt^\beta) &\leq \kappa (\vt) \leq \Ov{\kappa}(1 + \vt^\beta), 
	\quad \mbox{ where }\ \beta > 6. \label{w16}
\end{align}

As a consequence of the hypotheses \eqref{w10}, \eqref{w11}, \eqref{w13}, and \eqref{w14}, we get the following estimates, for which we refer
to \cite[Chapter 3, Section 3.2]{FeNo6A}:
\begin{align} 
	\vr^{\frac{5}{3}} + \vt^4 \lesssim \vr e(\vr, \vt) &\lesssim 	1+ \vr^{\frac{5}{3}} + \vt^4, \label{L5b} \\
	s_{\rm m}(\vr, \vt) &\lesssim \left( 1 + |\log(\vr)| + [\log(\vt)]^+ \right). \label{L5a}
\end{align} 

\subsection{Coercivity of the dissipative stress}

In order to obtain uniform energy bounds, we need to control the $L^2$-norm of $\vu$ in terms of the dissipation potential 
\begin{align}
	\intO{ \frac{1}{\vt} \mathbb{S}(\vt, \Ds \vu): \Ds \vu } &\approx \intO{\frac{\mu(\vt)}{\vt} |\Ds \vu - \frac{1}{3} \Div \vu \mathbb{I} |^2} +  \intO{\frac{\eta(\vt)}{\vt} |\Div \vu |^2} 
	\br &\gtrsim \underline{\mu} \intO{|\Ds \vu - \frac{1}{3} \Div \vu \mathbb{I} |^2} + \underline{\eta} \intO{ |\Div \vu|^2}. 
	\label{est:coerc}
\end{align}
This can be done only far away from the set
\[
\mathcal{N} = \left\{ \vc{w} \in W^{1,2}(\Omega; R^3) \Big| \ \underline{\mu} \left( \Ds \vc{w} - \frac{1}{3} \Div \vc{w} \mathbb{I}  \right) = 0,\quad
\underline{\eta} \Div \vc{w} = 0 \ \mbox{ in }\ \Omega, \qquad \vc{w} \cdot \vc{n}|_{\partial \Omega} = 0 \right\}.
\]
As a matter of fact, let $\Pi_{\mathcal{N}}$ denote a projection onto the space $\mathcal{N}$. Then
Korn--Poincar\' e inequality, see \tsl{e.g.} Lewintan, M\" uller, Neff \cite{LewMulNef}, gives
\begin{equation} \label{KoPo}
\| \vc{v} - \Pi_{\mathcal{N}} \vc{v} \|_{W^{1,2}(\Omega; R^3)} \aleq \left( \underline{\mu}\left\| \Ds \vc{v} - \frac{1}{3} \Div \vc{v} \mathbb{I} \right\|_{L^2(\Omega; R^{3 \times 3})} + 
\underline{\eta} \| \Div \vc{v} \|_{L^2(\Omega)} \right)
\end{equation}
for any $\vc{v} \in W^{1,2}(\Omega;R^3)$, $\vc{v} \cdot \vc{n}|_{\partial \Omega} = 0$. 

Notice that, from Schirra \cite[Proposition 2.5]{Schi}, it follows that, if $\underline{\mu} > 0$ and $\underline{\eta} = 0$, the set $\mathcal{N}$ contains all conformal Killing vectors with vanishing normal trace.
If 
instead $\underline{\mu}>0$ and
$\underline{\eta} > 0$, then the set $\mathcal{N}$ consists of all \emph{rigid motions} $\vc{w} = \bfomega \times x + \vc{a}$ with zero normal trace. In this
latter case, $\mathcal{N} \ne \{ 0 \}$ if and only if $\Omega$ is rotationally symmetric.

Finally, we point out that,
if the velocity field $\vu$ satisfies the complete slip boundary conditions \eqref{i5}, uniform energy estimates can be obtained provided the harmonic extension 
$\vtB$ of the boundary condition \eqref{i6} inside $\Omega$ satisfies 
\begin{equation} \label{COC} 
	\Grad \vtB \cdot \vc{w} = 0 \qquad \mbox{ for any }\quad \vc{w} \in \mathcal{N}.
\end{equation}
We refer to Section \ref{EE} below for more details.
Note that \eqref{COC} is satisfied, for instance, in the case when $\Omega$ is an annulus bounded by two concentric spheres, 
\[
\Omega = \left\{ x \in R^3 \ \Big|\ 0 < r_1 < |x| < r_2 \right\}\qquad \mbox{ and } \qquad \vtB = \Theta_1 \ \mbox{ if } \ |x| = r_1,\quad
\vtB = \Theta_2 \ \mbox{ if }\ |x| = r_2,
\]
provided $\underline{\mu}, \underline{\eta} > 0$, and $\Theta_1$, $\Theta_2$ are positive constants. 

\subsection{Existence of global in time weak solutions}

Under the hypotheses \eqref{w9}--\eqref{w16} and \eqref{COC}, the NSF system \eqref{i1}--\eqref{i8}, together with the boundary conditions
\eqref{i5}, \eqref{i6}, admits a global in time weak solution in the sense of Definition \ref{DL1} 
for any finite energy initial datum and any sufficiently smooth boundary datum $\vtB$, see \cite[Chapter 12, Theorem 18]{FeiNovOpen}.

\begin{Remark} \label{Rm1}
	
	Strictly speaking, the existence theory in \cite{FeiNovOpen} is formulated for the Dirichlet boundary conditions 
	for the velocity. The extension to the complete slip conditions \eqref{i5} is straightforward. 
	
	\end{Remark}
	
In addition, the weak solutions in the sense of Definition \ref{DL1} enjoy the \emph{weak--strong uniqueness} 
property; they coincide with the strong solution associated to the same initial/boundary data as long as the 
latter exists, see \cite[Section 4]{FeGwKwSG}.

\section{Main result: the singular limit problem}
\label{m}

Having collected all preliminary material, we are able to state the main result of this paper.

\begin{mdframed}[style=MyFrame]

\begin{Theorem} [{\bf Singular limit}] \label{TM1}
	
Let $\Omega \subset R^3$ be a bounded domain of class $C^3$. Let the constitutive relations for $p$, $e$, $s$, $\mu$, $\lambda$, 
and $\kappa$ comply with hypotheses \eqref{con1}--\eqref{w16}. Let $(\vre, \vue, \vte)_{\ep > 0}$ be a family of weak solutions 
to the Navier--Stokes--Fourier system, in the sense specified in Definition \ref{DL1}, emanating from the initial data 
\begin{align}
\vre(0, \cdot) = \Ov{\vr} + \ep \vr_{0,\ep},\ \intO{ \vr_{0,\ep} } = 0,\ \Ov{\vr} > 0, \ \vr_{0,\ep} \to 
\mathfrak{R}_0 \ \mbox{weakly-(*) in}\ L^\infty (\Omega), \br 
\vte(0, \cdot) = \Ov{\vt} + \ep \vt_{0,\ep},\ \Ov{\vt} > 0, \intO{ \vt_{0,\ep} } = 0,\ \vt_{0,\ep} \to 
	\mathfrak{T}_0 \ \mbox{weakly-(*) in}\ L^\infty (\Omega), \br 
\vue(0, \ep) = \vu_{0,\ep} \to \vu_0 \ \mbox{weakly-(*) in}\ L^\infty(\Omega; R^3).	 
\label{dt1}	
	\end{align}
In addition, the Dirichlet boundary conditions for the temperature are given, 
\begin{equation} \label{dt2}
\vte|_{\partial \Omega} = \Ov{\vt} + \ep \mathfrak{T}_B, 
\end{equation}
where the harmonic extension of $\mathfrak{T}_B$ satisfies the coercivity hypothesis \eqref{COC}. 

Then there is a subsequence (not relabelled) such that:

\begin{align} 
	\frac{\vre - \Ov{\vr}}{\ep} &\to \mathfrak{R} \ \mbox{weakly-(*) in} \ L^\infty(0,T; L^{\frac 5 3}(\Omega)), \br
	\frac{\vte - \Ov{\vt}}{\ep} &\to \mathfrak{T} \ \mbox{weakly-(*) in} \ L^\infty(0,T; L^{2}(\Omega)), \ \mbox{and weakly in}\ L^2(0,T; W^{1,2}(\Omega)), \br
	\vue &\to \vU \ \mbox{weakly in}\ L^2((0,T) \times \Omega; R^3),
\label{dt3}
\end{align}
where $\mathfrak{R}$, $\vU$, and 
\[
\Theta := \mathfrak{T} - \frac{\lambda(\Ov{\vr}, \Ov{\vt})}{|\Omega|} \intO{ \mathfrak{T}}
\]
represent a weak solution of the Oberbeck--Boussinesq system \eqref{ObBs1}--\eqref{i14B}, emanating from the initial data 
\begin{align}
\vU(0, \cdot) &= \vc{H}[\vu_0], \br
\Ov{\vr} c_p (\Ov{\vr}, \Ov{\vt}) \Theta(0, \cdot) &= 
\Ov{\vt} \left( \frac{\partial s (\Ov{\vr}, \Ov{\vt})}{\partial \vr} \mathfrak{R}_0 +  \frac{\partial s (\Ov{\vr}, \Ov{\vt})}{\partial \vt} \mathfrak{T}_0 + 
\alpha (\Ov{\vr}, \Ov{\vt}) G \right),
\label{dt4}
\end{align}
where $\vc{H}$ denotes the Helmholtz projection onto the space of solenoidal functions.

\end{Theorem}

\end{mdframed}

The rest of the paper is devoted to the proof of Theorem \ref{TM1}. As already mentioned, the next section is devoted to the derivation of suitable
uniform bounds for the family of weak solutions $(\vre, \vte, \vue)_{\ep > 0}$ to the scaled NSF system. In Section \ref{w}, we establish weak convergence
properties and compute the limit of most of the terms appearing in equations \eqref{i1}-\eqref{i3}. Finally, in Section \ref{A} we
compute the limit of the convective term appearing in the momentum equation \eqref{i2}, thus completing the proof to Theorem \ref{TM1}.

\section{Uniform bounds}
\label{e}

Given a global--in--time weak solution $(\vre, \vte, \vue)_{\ep > 0}$ to the scaled NSF system, 
our first goal is to establish uniform bounds, which are independent of the scaling parameter $\ep$. 

\subsection{Mass conservation}

As $\vu \cdot \vc{n}|_{\partial \Omega} = 0$, the total mass of the fluid is a conserved quantity. Consequently, in accordance with \eqref{S2}, \eqref{S3}, one has
\begin{equation} \label{ee1}
	\intO{ \vre (t,\cdot)} = \intO{ \vr_{0,\ep}} = \Ov{\vr}|\Omega|.
\end{equation}

\subsection{Essential vs. residual component}

Following 
\cite{FeNo6A}, we consider the \emph{essential} and \emph{residual} components of a measurable function $g=g(t,x)$. For a compact set 
\[
	K \subset \left\{ (\vr, \vt) \in \R^2 \ \Big| \ \vr > 0, \vt > 0 \right\}
\]
and $\ep > 0$, we introduce the functions
\[
	[g]_{\rm ess} = g \mathds{1}_{(\vre, \vte) \in K},\ 
	[g]_{\rm res} = g - [g]_{\rm ess} = g \mathds{1}_{(\vre, \vte) \in \R^2 \setminus K},
\]
where $\mathds{1}_A$ denotes the characteristic function of a set $A\subset\R^2$.

More specifically, we fix $\mathcal{U} (\Ov{\vr}, \Ov{\vt} )\subset (0, \infty)^2$ to be an open neighborhood of $(\Ov{\vr}, \Ov{\vt})$.
Note that $\mathcal{U}$ contains the range of $(\tvr_\ep, \Ov{\vt} + \ep \mathfrak{T}_B)$ for all $\ep > 0$ small enough.
Then, the set $K$ is chosen so that
\begin{equation} \label{ee2}
K = \Ov{\mathcal{U}(\Ov{\vr}, \Ov{\vt})}. 
\end{equation}

Next, we record the following bounds shown in \cite[Chapter 5, Lemma~5.1]{FeNo6A}:
\begin{align} 
	\left[E_{\ep} \left( \vre, \vte, \vue \Big| \tvr, \tvt, \tvu \right)\right]_{\rm ess} &\geq 
	C \left[ \frac{ |\vre - \tvr|^2 }{\ep^2} + \frac{ |\vte - \tvt|^2 }{\ep^2} + |\vue - \tvu |^2 \right]_{\rm ess}, 
	\label{BB1} \\ 
	\left[ E_{\ep} \left( \vre, \vte, \vue \Big| \tvr, \tvt, \tvu \right) \right]_{\rm res} &\geq 
	C  \left[ \frac{1}{\ep^2} + \frac{1}{\ep^2} \vre e(\vre, \vte) + \frac{1}{\ep^2} \vre |s(\vre, \vte)| + \vre |\vue|^2 \right]_{\rm res}, 
	\label{ee3}	
\end{align}
whenever $K$ is given by \eqref{ee2} and $(\tvr, \tvt) \in \mathcal{U}(\Ov{\vr}, \Ov{\vt})$. The constant $C$ depends on $K$ and on the distance
\[
\sup_{t,x} {\rm dist} \left[ (\tvr (t,x), \tvt (t,x) ) ; \partial K \right]. 	
\]

\subsection{Energy estimates}
\label{EE}

Using the ansatz $\tvr = \tvr_\ep$, $\tvt = \Ov{\vt} + \ep \mathfrak{T}_B$, $\tvu = 0$ in the relative energy inequality \eqref{L4}, we obtain
 \begin{align}
	&\left[ \intO{ E_\ep \left(\vre, \vte, \vue \Big| \tvr_\ep, \Ov{\vt} + \ep \mathfrak{T}_B, 0 \right) } \right]_{t = 0}^{t = \tau} + \frac{1}{\ep^2} \int_0^\tau \int_{\Ov{\Omega}} (\Ov{\vt} + \ep \mathfrak{T}_B)  \D \sigma_\ep (t,x) \br 
	&\leq - \frac{1}{\ep} \int_0^\tau \intOe{ \left(  \vre \Big(s(\vre, \vte) - s(\tvr_\ep, \Ov{\vt} + \ep \mathfrak{T}_B) \Big) \vue \cdot \Grad \mathfrak{T}_B -
		\left( \frac{\kappa (\vte) \Grad \vte}{\vte} \right) \cdot \Grad \mathfrak{T}_B \right) } \dt \br 
	&+ \int_0^\tau \intO{ \frac{1}{\ep} \vre \Grad G  \cdot \vue  } \dt  - \frac{1}{\ep^2} \int_0^\tau \intO{  \frac{\vre}{\tvr_\ep} \vue \cdot \Grad p(\tvr_\ep, \Ov{\vt} + \ep \mathfrak{T}_B)  } \dt,
	\label{e1}
\end{align}
where, according to \eqref{Lw7a}, one has
\begin{equation} \label{e2}
\sigma_\ep \geq 	\frac{1}{\vte} \left(  \ep^2 \mathbb{S}(\vte, \Ds \vue) : \Ds \vue + \frac{\kappa (\vte) |\Grad \vte |^2}{\vte}\right).	
\end{equation}	
Observe that, owing to \eqref{S1} and \eqref{T_B-harm}, we get
\begin{equation} \label{e3}
	\int_0^\tau \intO{ \frac{1}{\ep}  \vre \Grad G  \cdot \vue  } \dt - \frac{1}{\ep^2} \int_0^\tau \intOe{  \frac{\vre}{\tvr_\ep} \vue \cdot \Grad p(\tvr_\ep, \Ov{\vt} + \ep \mathfrak{T}_B)  } \dt = 0.
	\end{equation}
Thus, inequality \eqref{e1} reduces to 
 \begin{align}
	&\left[ \intO{ E_\ep \left(\vre, \vte, \vue \Big| \tvr_\ep, \Ov{\vt} + \ep \mathfrak{T}_B, 0 \right) } \right]_{t = 0}^{t = \tau} + \frac{1}{\ep^2} \int_0^\tau \int_{\Ov{\Omega}} (\Ov{\vt} + \ep \mathfrak{T}_B)  \D \sigma_\ep (t,x) \br 
	&\leq - \frac{1}{\ep} \int_0^\tau \intOe{ \left(  \vre \Big(s(\vre, \vte) - s(\tvr_\ep, \Ov{\vt} + \ep \mathfrak{T}_B) \Big) \vue \cdot \Grad \mathfrak{T}_B -
		\left( \frac{\kappa (\vte) \Grad \vte}{\vte} \right) \cdot \Grad \mathfrak{T}_B \right) } \dt. 
	\label{e4}
\end{align}
We are now going to bound each term appearing in the right--hand side of the previous relation.

\subsubsection{Entropy-dependent term}

Let us start by considering the entropy-dependent term. We resort to the decomposition into essential and residual sets to write
\begin{align}
	\frac{1}{\ep} &\intO{ \left|  \vre (s(\vre, \vte) - s(\tvr_\ep, \Ov{\vt} + \ep \mathfrak{T}_B)) \vue \cdot \Grad \mathfrak{T}_B \right| }	\br 
	&\lesssim \frac{1}{\ep} \intO{ \left| \left[ \vre (s(\vre, \vte) - s(\tvr_\ep, \Ov{\vt} + \ep \mathfrak{T}_B))\vue \right]_{\rm ess}  \right| } \br &+ 
	\frac{1}{\ep} \intO{ \left| \left[ \vre (s(\vre, \vte) - s(\tvr_\ep, \Ov{\vt} + \ep \mathfrak{T}_B)) \vue \cdot \Grad \mathfrak{T}_B  \right]_{\rm res}  \right| }.
	\nonumber
\end{align}
As for the essential part term, we have
\begin{align}
	\frac{1}{\ep} &\intO{ \left| \left[ \vre (s(\vre, \vte) - s(\tvr_\ep, \Ov{\vt} + \ep \mathfrak{T}_B)) \vue \right]_{\rm ess}  \right| } \br
&\lesssim
	\frac{1}{\ep^2} \intO{ \left| \left[ (s(\vre, \vte) - s(\tvr_\ep, \Ov{\vt} + \ep \mathfrak{T}_B)) \right]_{\rm ess} \right|^2  } + \intO{ \vre |\vue|^2 } \br
&\lesssim \intO{ E_\ep \left( \vre, \vte, \vue \Big| \tvr_\ep, \Ov{\vt} + \ep \mathfrak{T}_B, 0 \right) }
	\label{e5}
\end{align}
by virtue of \eqref{S4} and \eqref{BB1}, while the residual part term can be controlled as
\begin{align} 
	\frac{1}{\ep} &\intO{ \left| \left[ \vre (s(\vre, \vte) - s(\tvr_\ep, \Ov{\vt} + \ep \mathfrak{T}_B)) \vue \cdot \Grad \mathfrak{T}_B \right]_{\rm res}  \right| }\br
	&\lesssim
	\frac{1}{\ep} \intO{ \left[ \vre  |\vue| \right]_{\rm res} } + 
	\frac{1}{\ep} \intO{ \left[ \vre s_{\rm m}(\vre, \vte) |\vue| \right]_{\rm res} } + 
	\frac{1}{\ep} \intO{ \left[  \vte^3 |\vue \cdot  \Grad \mathfrak{T}_B | \right]_{\rm res} }. 
	\nonumber	
\end{align}

In view of \eqref{ee1}, the total mass is conserved and we may infer that 
\begin{equation} \label{e6}
	\frac{1}{\ep} \intO{ \left[ \vre  |\vue| \right]_{\rm res} } \lesssim \frac{1}{\ep^2} \intO{ [\vr_\ep]_{\rm res} } + \intO{\vre |\vue|^2} \lesssim
	\intO{ E_\ep \left( \vre, \vte, \vue \Big| \tvr_\ep, \Ov{\vt} + \ep \mathfrak{T}_B, 0 \right) }.
\end{equation}
Furthermore, by virtue of \eqref{L5b} and \eqref{L5a}, we have
\begin{align} 
	\frac{1}{\ep} \intO{ \left[ \vre s_{\rm m}(\vre, \vte) |\vue| \right]_{\rm res} } &\lesssim \frac{1}{\ep^2} \intO{ \left[ \vre s^2_{\rm m}(\vre, \vte) \right]_{\rm res} } + \intO{\vre |\vue|^2}  \br
	&\lesssim
	\intO{ E_\ep \left( \vre, \vte, \vue \Big| \tvr_\ep, \Ov{\vt} + \ep \mathfrak{T}_B, 0 \right) }.
	\label{e7}	
\end{align}

Finally, we introduce the following refinement of the residual component of a function: for a given $\theta_*\gg1$ to be fixed later, we define
\[
[g]_{\rm res,S} = g \mathds{1}_{(\vre, \vte) \in \R^2 \setminus K, \vte\leq \theta_*}, \qquad
[g]_{\rm res,L} = g \mathds{1}_{(\vre, \vte) \in \R^2 \setminus K, \vte\geq \theta_*},
\]
where the index $S$ stay for ``small'' and the index $L$ for ``large''. Using that
\[
[g]_{\rm res} = [g]_{\rm res,S} + [g]_{\rm res,L},
\]
by simple computations and Cauchy--Schwarz inequality we get 
\begin{align} 
	\frac{1}{\ep} \intO{ \left[  \vte^3 |\vue \cdot \Grad \mathfrak{T}_B| \right]_{\rm res} } &\lesssim 2\delta \intO{ |\vue \cdot \Grad \mathfrak{T}_B |^2} \br
&\quad + \frac{C(\delta,\theta_*)}{\veps^2}\intO{[1]_{\rm res}}	+ \frac{  C(\delta) }{\ep^2} \intO{ [\vte^6]_{\rm res,L} }  \br
&\lesssim 2\delta \intO{ |\vue \cdot \Grad \mathfrak{T}_B |^2} \br
&\quad + C(\delta,\theta_*) \intO{ E_\ep \left( \vre, \vte, \vue \Big| \tvr_\ep, \Ov{\vt} + \ep \mathfrak{T}_B, 0 \right) }	\br 
&\qquad\qquad+ 
\frac{  C(\delta) }{\ep^2} \intO{ [\vte^6]_{\rm res,L} }
	\label{e8}
\end{align}
for any $\delta > 0$. Observe that we have also used \eqref{ee3} for passing from the first to the second inequality.
First, we rewrite 
\[
\intO{ |\vue \cdot \Grad \mathfrak{T}_B |^2 } = 
\intO{ |(\vue - \Pi_{\mathcal{N}} \vue ) \cdot \Grad \mathfrak{T}_B |^2 } + 
\intO{ | \Pi_{\mathcal{N}} \vue\cdot \Grad \mathfrak{T}_B |^2 }, 
\]
where, in accordance with hypothesis \eqref{COC}, one has
\[
\Pi_{\mathcal{N}} \vue\cdot \Grad \mathfrak{T}_B = 0.
\]
Second, by virtue of Korn--Poincar\' e inequality \eqref{KoPo} and \eqref{est:coerc}, we get 
\[
\intO{ |(\vue - \Pi_{\mathcal{N}} \vue ) \cdot \Grad \mathfrak{T}_B |^2 } \lesssim \intO{ \frac{1}{\vte} \mathbb{S}(\vte, \Ds \vue): \Ds \vue }. 
\]
Consequently, the first integral on the right--hand side of \eqref{e8} can be absorbed by the left--hand side of \eqref{e4} as soon as $\delta > 0$ is fixed small enough.

Let us now deal with the last integral in \eqref{e8}. 
We start by writing 
the bound
\begin{equation} \label{est:theta^6}
\frac{ 1 }{\ep^2} \intO{ [\vte^6]_{\rm res,L} } \lesssim
\frac{1}{\ep^2} \intO{ \left( [\vte^3 - \theta_*^3 ]^+ \right)^2 } + \frac{1}{\ep^2} \intO{ [\theta_*^6]_{\rm res} },
\end{equation}
where, by virtue of \eqref{ee3}, 
one has 
\begin{equation} \label{est:theta-*}
\frac{1}{\ep^2} \intO{ [\theta_*^6]_{\rm res} } \lesssim \intO{ E_\ep \left( \vre, \vte, \vue \Big| \tvr_\ep, \Ov{\vt} + \ep \mathfrak{T}
	_B, 0 \right) }.
\end{equation}
Next, 
by virtue of Poincar\' e inequality, one has
\[
\frac{1}{\ep^2} \intO{ \left( [\vte^3 - \theta_*^3 ]^+ \right)^2 } 
\lesssim \frac{1}{\ep^2} \int_{\{\vte\geq \theta_*\}}{|\Grad \vte^3|^2 }\dx 
\]
where no boundary term appears on the right-hand side provided we choose (say) $\theta_*\geq 2\Ov{\vt}$.
At this point, as $\beta > 6$, we can bound
\begin{align*}
\frac{1}{\ep^2} \int_{\{\vte\geq \theta_*\}}{|\Grad \vte^3|^2 }\dx &\approx \frac{1}{\ep^2} \int_{\{\vte\geq \theta_*\}}{ \vte^4 |\Grad \vte|^2 }\dx \\
&\leq \frac{\theta_*^{-(\beta-6)}}{\ep^2} \intO{ \frac{\kappa (\vte) |\Grad \vte|^2 }{\vte^2} }.
\end{align*}
In particular, inserting all these bounds into \eqref{est:theta^6} and then into \eqref{e8} and finally choosing $\theta_*\gg1$ large enough,
owing to \eqref{e2} the integral on the right--hand side of the previous
relation can be absorbed by the left--hand side of \eqref{e4}. 

Summing up the previous estimates, we may rewrite inequality \eqref{e4} in the form 
 \begin{align}
	&\left[ \intO{ E_\ep \left(\vre, \vte, \vue \Big| \tvr_\ep, \Ov{\vt} + \ep \mathfrak{T}_B, 0 \right) } \right]_{t = 0}^{t = \tau} + \frac{1}{\ep^2} \int_0^\tau \int_{\Ov{\Omega}} (\Ov{\vt} + \ep \mathfrak{T}_B)  \D \sigma_\ep (t,x) \br 
	&\lesssim \left( 1 + \frac{1}{\ep} \left| \int_0^\tau \intO{ \frac{\kappa (\vte) \Grad \vte}{\vte} \cdot \Grad\mathfrak{T}_B  } \dt \right| + 
	\int_0^\tau \intO{ E_\ep \left(\vre, \vte, \vue \Big| \tvr_\ep, \Ov{\vt} + \ep \mathfrak{T}_B, 0 \right) } \dt \right).
 	\label{e9}
\end{align}

\subsubsection{Heat flux dependent term}

We are now going to bound the second term appearing in the right--hand side of \eqref{e9}.
Evoking once more hypothesis \eqref{w16} we get 
\[
\frac{1}{\ep} \left| \intO{ \frac{\kappa(\vte) }{\vte} \Grad \vte \cdot \Grad\mathfrak{T}_B } \right| \lesssim
\frac{1}{\ep} \intO{ |\Grad (\log(\vte))| + \vte^{\beta - 1} |\Grad \vte| },
\]
where we can bound
\begin{align}
	\frac{1}{\ep} \intO{ \left|  \Grad (\log(\vte)) \right| } \lesssim \frac{\delta}{\ep^2} \intO{ \left|  \Grad (\log(\vte)) \right|^2 } + C(\delta)
	\nonumber
\end{align}
for any $\delta> 0$ to be fixed later. Thus, the integral on the right--hand side is controlled by the left--hand side of \eqref{e9},
provide we choose $\de>0$ small enough.

Next, we compute
\begin{align}
	\frac{1}{\ep} \intO{ \vte^{\beta - 1}{\Grad \vte} } &= 
	\frac{1}{\ep} \intO{ \vt^{\frac{\beta}{2}} \vte^{\frac{\beta}{2}-1} \Grad \vte }\leq \frac{\delta}{\ep^2} \intO{ \frac{\kappa (\vte) |\Grad \vte|^2 }{\vte^2} } + C(\delta) \intO{ |\vte^{\frac{\beta}{2} }|^2 }, 	
	\nonumber
\end{align}
where the first term is absorbed by the left--hand side of \eqref{e9}. Finally, using Poincar\' e inequality as done above, we infer 
\begin{align}
	\intO{ |\vte^{\frac{\beta}{2} } |^2 } \lesssim \intO{ |\Grad \vte^{\frac{\beta}{2}} |^2 } + \int_{\partial \Omega} ( \Ov{\vt} + \ep \mathfrak{T}_B )^\beta \ \D S,
	\nonumber
\end{align}
where one has
\[
\intO{ |\Grad \vte^{\frac{\beta}{2}} |^2 } = \intO{ \vte^{\bt-2} |\Grad \vte |^2 } \leq
\intO{ \frac{\kappa (\vte) |\Grad \vte|^2 }{\vte^2} }.
\]

In the end, we may apply Gr\"onwall's lemma to \eqref{e9} to conclude that, for any time $T>0$ fixed, there holds
\begin{equation} \label{e10}
\sup_{t \in [0,T]}  \intO{ E_\ep \left(\vre, \vte, \vue \Big| \tvr_\ep, \Ov{\vt} + \ep \mathfrak{T}_B, 0 \right) } + 
\frac{1}{\ep^2} \int_0^T \int_{\Ov{\Omega}} (\Ov{\vt} + \ep \mathfrak{T}_B)  \D \sigma_\ep (t,x) \lesssim 1 
\end{equation}
as long as 
\begin{equation} \label{e11}
 \intO{ E_\ep \left(\vre, \vte, \vue \Big| \tvr_\ep, \Ov{\vt} + \ep \mathfrak{T}_B, 0 \right) (0, \cdot) } \lesssim 1, 
\end{equation}
meaning the initial data are ill--prepared, cf. hypotheses \eqref{wpd}, \eqref{wpd1}. We point out that the (implicit) multiplicative
constant in \eqref{e10} depends on the fixed time $T>0$.

\subsection{Uniform bounds -- summary}

The bounds established in \eqref{e10}, \eqref{e11}, together with the structural hypotheses \eqref{w9}--\eqref{w16}, yield the following uniform estimates
for $\ep \to 0$, see \cite[Chapter 5, Proposition 5.1]{FeNo6A}:

\begin{align}
	{\rm ess} \sup_{t \in (0,T)} \intO{ \mathds{1}_{\rm res}(t, \cdot)   } &\lesssim \ep^2 \label{ub1}, \\
	{\rm ess} \sup_{t \in (0,T)} \left\| \left[\frac{\vre - \Ov{\vr}}{\ep} \right]_{\rm ess} (t, \cdot) \right\|_{L^2(\Omega)} &\lesssim 1, \label{ub2}\\
		{\rm ess} \sup_{t \in (0,T)} \left\| \left[\frac{\vte - \Ov{\vt}}{\ep} \right]_{\rm ess} (t, \cdot) \right\|_{L^2(\Omega)} &\lesssim 1, \label{ub3}\\
			{\rm ess} \sup_{t \in (0,T)} \intO{ \left( [\vr]_{\rm res}^{\frac{5}{3}} + [\vt]_{\rm res}^{4}  \right)(t, \cdot) } &\lesssim \ep^2, \label{ub4} \\ 
			{\rm ess} \sup_{t \in (0,T)} \left\| \sqrt{\vre} \vue (t, \cdot) \right\|_{L^2(\Omega;R^3)}	&\lesssim 1, \label{ub5} \\ 
			\int_0^T \int_{\Ov{\Omega}} 1 \ \D \sigma (t,x) &\lesssim \ep^2, \label{ub5a} \\ 
			\int_0^T \| \vue (t, \cdot) \|^2_{W^{1,2}(\Omega; R^3)} \dt &\lesssim 1, \label{ub6}\\
			\int_0^T \left( \left\| \frac{\vte - \Ov{\vt} }{\ep} (t, \cdot) \right\|_{W^{1,2}(\Omega)}  +  \left\| \frac{\log(\vte) - \log(\Ov{\vt}) }{\ep} (t, \cdot) \right\|_{W^{1,2}(\Omega)}    \right) &\lesssim 1, \label{ub7} \\ 
			\int_0^T \left\| \left[ \frac{\vre s(\vre, \vte)}{\ep}    \right]_{\rm res} (t, \cdot) \right\|_{L^q(\Omega)}^q \dt &\lesssim 1 
			\ \mbox{for a certain}\ q > 1, \label{ub8}\\
				\int_0^T \left\| \left[ \frac{\vre s(\vre, \vte)}{\ep}    \right]_{\rm res} \vue (t, \cdot) \right\|_{L^q(\Omega;R^3)}^q \dt &\lesssim 1 
			\ \mbox{for a certain}\ q > 1, \label{ub9}\\
			\int_0^T \left\| \left[ \frac{\kappa (\vte)}{\vte}\right]_{\rm res} \left( \frac{\Grad \vte}{\ep }\right) (t, \cdot) \right\|_{L^q(\Omega;R^3)}^q \dt
			&\lesssim 1 . 
			\label{ub10}
		\end{align}

\section{Weak convergence towards the target system}
\label{w}

Making use of the uniform bounds \eqref{ub1}--\eqref{ub10} we perform the limit in the weak topologies. This part of the proof, with the exception of the limit in the momentum equation, is almost 
identical with \cite[Chapter 5, Section 5.3]{FeNo6A}. We therefore only state the results referring to \cite{FeNo6A} for details. 

First, it follows from the uniform bounds \eqref{ub2}, \eqref{ub4} that 
\begin{equation} \label{w1}
	\frac{\vre - \Ov{\vr}}{\ep} \to \mathfrak{R} \qquad \mbox{ weakly-(*) in } \quad L^\infty(0,T; L^{\frac{5}{3} }(\Omega)),  
	\end{equation}
in particular we have the strong convergence
\begin{equation}\label{w2}
	\vre \to \Ov{\vr}\qquad \mbox{ in }\quad 	L^\infty(0,T; L^{\frac{5}{3} }(\Omega)).
\end{equation}

Similarly, we gather
\begin{align} 
	\frac{\vte - \Ov{\vt}}{\ep} &\to \mathfrak{T} \qquad \mbox{ weakly-(*) in } \quad L^\infty(0,T; L^{2}(\Omega)), \quad \mbox{ weakly in }\quad L^2(0,T; W^{1,2}(\Omega)), \label{w3}\\
	\vte &\to \Ov{\vt} \qquad \mbox{ in }\quad L^\infty(0,T; L^{2}(\Omega)). \label{w4}  
\end{align}

Finally, by virtue of \eqref{ub6}, we get
\begin{equation} \label{w5}
	\vue \to \vU \qquad \mbox{ weakly in }\quad L^2(0,T; W^{1,2}(\Omega;R^3) ).
	\end{equation}
In all cases, we have to extract a suitable subsequence as the case may be.

In accordance with the boundary conditions \eqref{i5}, \eqref{i6}, we have
\begin{align} 
	\mathfrak{T}|_{\partial \Omega} &= \mathfrak{T}_B , \label{w6}\\
\vU \cdot \vc{n}|_{\partial \Omega} &=0 . \label{w7}
\end{align}
Furthermore, it follows from \eqref{ee1} that 
\begin{equation} \label{w8}
	\intO{ \mathfrak{R}(t, \cdot)} = 0 .
	\end{equation}

	\subsection{Equation of continuity}
	
With \eqref{w2}, \eqref{w5} at hand, we may let $\ep \to 0$ in the weak formulation of the equation of continuity \eqref{Lw4} to obtain 
\begin{equation} \label{ww9}
	\Div \vU = 0 \qquad \mbox{ a.a. in }\quad (0,T) \times \Omega.
\end{equation}

\subsection{Entropy equation}

Passing to the limit in the entropy equation \eqref{Lw7} is a bit lengthy but, fortunately, the same as in \cite[Chapter 5, Section 5.3.2]{FeNo6A}. Thus we only report the result: 
\begin{align}
\int_0^T &\intO{\Ov{\vr} \left( \frac{\partial s(\Ov{\vr}, \Ov{\vt} ) }{\partial \vr } \mathfrak{R} + \frac{\partial s(\Ov{\vr}, \Ov{\vt} ) }{\partial \vt } \mathfrak{T}    \right)
\left( \partial_t \varphi  + \vU \cdot \Grad \varphi \right) } \dt	- \int_0^T \intO{ \frac{\kappa(\Ov{\vt})}{\Ov{\vt}} \Grad \mathfrak{T} \cdot \Grad \varphi} \dt \br
&= - \intO{\Ov{\vr} \left( \frac{\partial s(\Ov{\vr}, \Ov{\vt} ) }{\partial \vr } \mathfrak{R}_0 + \frac{\partial s(\Ov{\vr}, \Ov{\vt} ) }{\partial \vt } \mathfrak{T}_0    \right)\varphi (0, \cdot)  },
\label{ww10}
\end{align}
for any $\varphi \in C^1_c ([0, T) \times \Omega)$,
where we have set
\begin{equation} \label{ww11}
\vr_{0,\ep} \to \mathfrak{R}_0,\quad \vt_{0,\ep} \to \mathfrak{T}_0	\qquad \mbox{ weakly-(*) in }\quad L^\infty (\Omega).
	\end{equation}
Here, in order to conclude, we need a relation between the limit $\mathfrak{R}$ and $\mathfrak{T}$. This issue will be discussed in the forthcoming section.

\subsection{Momentum equation}

The asymptotic limit of the momentum equation \eqref{Lw5} is exactly the same as in \cite[Chapter 5, Section 5.5.3]{FeNo6A}, specifically, 
\begin{align}
\int_0^T &\intO{ \Big( \Ov{\vr} \vU \cdot \partial_t \bfphi + \Ov{\vr \vU \otimes \vU}: \Grad \bfphi - \mathbb{S}(\Ov{\vt}, \Ds \vU) : \Grad \bfphi + \mathfrak{R} \Grad G \cdot \bfphi \Big)  }\dt \br	
&= - \intO{ \Ov{\vr} \vU_0 \cdot \bfphi (0, \cdot)}
\label{ww12} 
 	\end{align}
for any test function 
\[
\bfphi \in C^1_c([0,T) \times \Ov{\Omega}; R^3),\quad \Div \bfphi = 0,\quad \bfphi \cdot \vc{n}|_{\partial \Omega} = 0.
\]
Here we have used the weak limits
\begin{align}
\vu_{0, \ep} &\to \vU_0 \qquad \mbox{ weakly in } \quad L^\infty(\Omega; R^3), \br 
\vre \vue \otimes \vue &\to \Ov{\vr \vU \otimes \vU} \qquad \mbox{ weakly in }\quad L^2(0,T; L^q(\Omega;R^3)),\ \mbox{for a certain}\ q > 1.
\label{ww13}
\end{align}

\subsubsection{Pressure term}

The limit equation \eqref{ww12} holds thanks to solenoidality of the test function. Multiplying the momentum equation by $\ep$ and performing the limit for a general 
test function $\bfphi$ we deduce, exactly as in \cite[Chapter 5, Section 5.5.3]{FeNo6A}, the Boussinesq relation 
\begin{equation} \label{ww14}
\frac{\partial p(\Ov{\vr}, \Ov{\vt})}{\partial \vr} \Grad \mathfrak{R} + \frac{\partial p(\Ov{\vr}, \Ov{\vt})}{\partial \vt} \Grad \mathfrak{T} = \Ov{\vr} \Grad G.
\end{equation}

Hereafter, we normalize the potential $G$ so that 
\begin{equation} \label{ww15}
\intO{G} = 0.
\end{equation}
Accordingly, using \eqref{w8} we may integrate \eqref{ww14} obtaining the desired relation between $\mathfrak{R}$ and $\mathfrak{T}$, namely 
\begin{equation} \label{ww16}
	\frac{\partial p(\Ov{\vr}, \Ov{\vt})}{\partial \vr} \mathfrak{R} + \frac{\partial p(\Ov{\vr}, \Ov{\vt})}{\partial \vt}  \mathfrak{T} = \Ov{\vr}  G + 
	\frac{\partial p(\Ov{\vr}, \Ov{\vt})}{\partial \vt} \frac{1}{|\Omega|} \intO{\mathfrak{T} }.
\end{equation}	

\subsubsection{Final form of the limit entropy equation}

The relation \eqref{ww16} plugged in the limit entropy equation \eqref{ww10} gives rise to the final form of the limit heat equation. 
Setting 
\begin{equation} \label{ww17}
\Theta = \mathfrak{T} - \frac{\lambda(\Ov{\vr}, \Ov{\vt})}{|\Omega|} \intO{ \mathfrak{T} } 
\end{equation}
we conclude that
\begin{align}
\Ov{\vr} c_p(\Ov{\vr}, \Ov{\vt} ) \Big( \partial_t \Theta + \vU \cdot \Grad \Theta \Big)	- 
\Ov{\vr} \ \Ov{\vt} \alpha(\Ov{\vr}, \Ov{\vt} ) \vU \cdot \Grad G
&= \kappa(\Ov{\vt}) \Del \Theta, \br 
\Theta|_{\partial \Omega} &=  \mathfrak{T}_B - \frac{\lambda(\Ov{\vr}, \Ov{\vt})}{1 - \lambda(\Ov{\vr}, \Ov{\vt})} \frac{1}{|\Omega|}\intO{ \Theta } .
 \label{ww18}
\end{align}

The initial value of $\Theta$ can be evaluated by means of \eqref{ww10}:
\[
\left( \frac{\partial s(\Ov{\vr}, \Ov{\vt} ) }{\partial \vr } \mathfrak{R} + \frac{\partial s(\Ov{\vr}, \Ov{\vt} ) }{\partial \vt } \mathfrak{T}    \right)(0, \cdot) = 
\left( \frac{\partial s(\Ov{\vr}, \Ov{\vt} ) }{\partial \vr } \mathfrak{R}_0 + \frac{\partial s(\Ov{\vr}, \Ov{\vt} ) }{\partial \vt } \mathfrak{T}_0    \right),
\]
where, by virtue of \eqref{ww16}, one has
\[
\mathfrak{R} (0, \cdot) = 	\left( \frac{\partial p(\Ov{\vr}, \Ov{\vt})}{\partial \vr} \right)^{-1} \left(\Ov{\vr} G + 	\frac{\partial p(\Ov{\vr}, \Ov{\vt})}{\partial \vt} \left( \frac{1}{|\Omega|} \intO{\mathfrak{T} (0, \cdot) } - \mathfrak{T}(0, \cdot)  \right) \right).
\]

To simplify, we suppose 
\begin{equation} \label{ww19}
	\intO{ \mathfrak{R}_0 } = \intO{ \mathfrak{T}_0 } = 0
	\end{equation}
yielding 
\[
\intO{ \mathfrak{T} (0,\cdot) = 0 }.
\]
Consequently, a bit lengthy but straightforward calculation yields
\begin{align}
\Ov{\vr} c_p (\Ov{\vr}, \Ov{\vt}) \Theta(0, \cdot) = 
\Ov{\vt} \left( \frac{\partial s (\Ov{\vr}, \Ov{\vt})}{\partial \vr} \mathfrak{R}_0 +  \frac{\partial s (\Ov{\vr}, \Ov{\vt})}{\partial \vt} \mathfrak{T}_0 -
\alpha (\Ov{\vr}, \Ov{\vt}) G \right).
\label{ww20}
\end{align}

\section{Propagation of acoustic waves}
\label{A}

To complete the proof of convergence, we have to identify the ``convective'' term $\Ov{\vr \vU \otimes \vU }$ in the momentum equation \eqref{ww12}.
Specifically, our goal is to show that
\begin{equation} \label{A1}
	\int_0^T \intO{ \Ov{\vr \vU \otimes \vU} : \Grad \bfphi} \dt = \int_0^T \intO{ \Ov{\vr} (\vU \otimes \vU): \Grad \bfphi } \dt
\end{equation}
for any test function
\begin{equation} \label{cond:test-f}
\bfphi \in C^1_c([0,T) \times \Ov{\Omega}; R^3),\quad \Div \bfphi = 0,\quad \bfphi \cdot \vc{n}|_{\partial \Omega} = 0.
\end{equation}

To begin with, let $\vc{H}$ be the Helmholtz projection onto the space of solenoidal functions with vanishing normal trace and consider
consider the Helmholtz decomposition 
\[
\vc{v} = \vc{H}[\vc v] + \vc{H}^\perp [\vc{v}] , 
\]
over $\R^3$.
As observed in \cite[Section 5.4.2]{FeNo6A}, the component
$\vc{H}[\vre \vue]$ is compact (in a suitable topology), in particular it converges almost everywhere on any set $(0,T)\times\Omega$.
Hence, the desired relation \eqref{A1} follows as soon as we show
\begin{equation} \label{A2}
	\int_0^T \intO{ \vc{H}^\perp [\vre \vue] \otimes \vc{H}^\perp [\vue] : \Grad \bfphi } \dt \to 0 \qquad \mbox{ as }\quad \ep \to 0
	\end{equation}
for any test function $\bfphi$ satisfying \eqref{cond:test-f}. 

\subsection{Acoustic equation}

The main problem in showing \eqref{A2} are rapid \emph{time oscillations} related to acoustic waves. The equation describing 
their evolution -- the acoustic equation -- has been derived in \cite[Chapter 5, Section 5.4.7]{FeNo6A}. Here, the main difficulty is the 
entropy equation holds only for all test functions which are compactly supported in $\Omega$, see \eqref{Lw7}.
Thus, our first goal is to
extend the validity of the entropy equation \eqref{i2} to test functions which may be non-zero at the boundary $\partial\Omega$. As we are going to see in a while,
the price to pay is the appearing of additional (small) terms in the weak formulation \eqref{Lw7}.

Let $\varphi (t,x) \in C^1_c ((0,T) \times \Ov{\Omega})$ be a given test function.
We consider its approximation 
\begin{equation} \label{A3}
\varphi_\ep(t,x) = \chi_\ep (x) \varphi(t,x),\qquad \chi_\ep \in \DC(\Omega),\ 0 \leq \chi_\ep \leq 1,\ 
\chi_\ep (x) = 1 \ \mbox{ whenever }\ {\rm dist}[x, \partial \Omega] > \ep.
\end{equation}
Obviously, $\varphi_\ep$ is a legitimate test function for the entropy balance  \eqref{Lw7}: using it in
that relation yields
	\begin{align}
	- \int_0^T &\intO{ \left[ \ep \vre \frac{ s(\vre, \vte) - s(\Ov{\vr}, \Ov{\vt}) }{\ep } \chi_\ep \partial_t \varphi + \ep \vre \frac{s (\vre ,\vte) -s(\Ov{\vr}, \Ov{\vt}) }{\ep } \vue \cdot \Grad \varphi + \frac{\vc{q} (\vte, \Grad \vte )}{\vte} \cdot 
		\Grad \varphi \right] } \dt \br &= \int_0^T \int_{\Omega}{ \varphi_\ep \ \D \sigma_\ep(t,x)} 
	+ \ep \int_0^T \intO{ \vre \frac{s (\vre ,\vte) -s(\Ov{\vr}, \Ov{\vt}) }{\ep } \vue \cdot \Grad (\varphi - \varphi_\ep)} \br
	&+ \int_0^T \intO{\frac{\vc{q} (\vte, \Grad \vte )}{\vte} \cdot 
	\Grad (\varphi - \varphi_\ep) }.	
	\label{A4} 
\end{align}
Consequently, the entropy balance \eqref{A4}, compared to its ``conservative'' counterpart in \cite[Chapter 5]{FeNo6A},
contains extra error terms represented by the last two integrals in \eqref{A4}. It is easy to check that
\begin{equation} \label{A5} 
	\Grad (\varphi - \varphi_\ep) = (1 - \chi_\ep) \Grad \varphi - \varphi \Grad \chi_\ep, 
	\end{equation}
where one has the bounds
\begin{equation} \label{A6}
\| \Grad \chi_\ep \|_{L^q(\Omega; R^3)} \lesssim \ep^{\frac{1 - q}{q}}, \ 1 \leq q \leq \infty	, \qquad
\| 1 - \chi_\ep \|_{L^p(\Omega)} \lesssim \ep^{\frac1 p} ,\ 1 \leq p \leq \infty.
\end{equation}

Using the error estimates \eqref{A6},
we may repeat step by step the arguments of \cite[Chapter 5, Section 5.4.7]{FeNo6A},
to deduce the following acoustic equation:

\begin{align} 
\int_0^T \intO{ \Big( \ep Z_\ep \partial_t \varphi + \vre \vue \cdot \Grad \varphi \Big) } \dt
&= \ep \int_0^T \intO{\vc{h}^1_\ep \cdot \Grad \varphi } + \ep^\g \int_0^T \intO{ h^2_\ep \varphi  } \dt \br
&\mbox{for any}\ \varphi \in C^1_c((0,T) \times \Ov{\Omega}),
\label{A7} \\
\int_0^T \intO{ \Big( \ep \vre \vue \cdot \partial_t \bfphi + \omega Z_\ep \Div \bfphi \Big) }  \dt &= 
\ep^\g \int_0^T \intO{\mathbb{H}^3_\ep : \Grad \bfphi   } \dt + \ep^\g \int_0^T \intO{\vc{h}^4_\ep \cdot \bfphi} \br
&\mbox{for any} \ \bfphi \in C^1_c((0,T) \times \Ov{\Omega}),\ \bfphi \cdot \vc{n}|_{\partial \Omega} = 0.
\label{A8}
\end{align}
Here we have $0<\g<1$ and we have defined
\begin{align} 
Z_\ep &= \frac{1}{\omega} \left(\omega \frac{\vre - \Ov{\vr}}{\ep} + A \vre \chi_\ep \frac{s(\vre, \vte) - s (\Ov{\vr}, \Ov{\vt})}{\ep} 
- \Ov{\vr} G + \frac{A}{\ep} \Sigma_\ep\right), \br
\Sigma_\ep (\tau, \cdot ) &= \int_0^\tau \chi_\ep \sigma_\ep (t, \cdot) \dt,
\nonumber
\end{align}
with the parameters $A$ and $\o$ defined by
\begin{align}
A &= \frac{1}{\Ov{\vr}} \frac{\partial p(\Ov{\vr}, \Ov{\vt})} {\partial \vt} \left(\frac{\partial s(\Ov{\vr}, \Ov{\vt})} {\partial \vt}
\right)^{-1} \quad \mbox{ and }\quad 
\omega &= \frac{\partial p(\Ov{\vr}, \Ov{\vt})} {\partial \vr} + 
\frac{1}{\Ov{\vr}^2} \left(\frac{\partial s(\Ov{\vr}, \Ov{\vt})} {\partial \vt} \right)^{-1} \left| \frac{\partial p(\Ov{\vr}, \Ov{\vt})} {\partial \vt} \right|^2.
\nonumber
\end{align}
In addition, we have the following bounds:
\begin{align}
\| \vc{h}^1_\ep	\|_{L^q(0,T; L^1(\Omega;R^3))} + \| {h}^2_\ep	\|_{L^q(0,T; L^1(\Omega))} + \| \vc{h}^4_\ep	\|_{L^q(0,T; L^1(\Omega;R^3))} & \aleq 1, \br
\| \mathbb{H}^3_\ep \|_{L^q(0,T; L^1(\Omega; R^{3 \times 3}))} & \aleq 1 
\ \mbox{for some}\ q > 1.
\label{A10}
	\end{align}

Notice that equations \eqref{A7}, \eqref{A8} correspond to the weak formulation of the following wave system:
\begin{equation} \label{eq:wave}
\left\{ \begin{array}{l}
         \ep \partial _tZ_\ep + \Div \vme = \ep \Div \vc{h}^1_\ep + \ep^\g h^2_\ep \\[1ex]
         \ep \partial_t\vme + \o \Grad Z_\ep = \ep^\g \Div \mathbb{H}^3_\ep + \ep^\g \vc{h}^4_\ep\,,
        \end{array}
\right.
\end{equation}
where we have set $\vme := \vre \vue$.

Observe that, in contrast with in \cite[Section 5.4.7]{FeNo6A}, the error terms represented 
by the right--hand side of \eqref{eq:wave} are larger, of order $\ep^\g$ with $\g\in(0,1)$, whereas $\g = 1$ in \cite[Section 5.4.7]{FeNo6A}.
Notice that a similar situation appeared in \cite{DS-F-S-WK, Fan}, for instance, in the context of the multiscale analysis
for geophysical flows.
We are going to apply a strategy similar to the one employed in those papers (see also \cite[Sections 5.4.5--5.4.7]{FeNo6A}),
based on Lions-Masmoudi compensated compactness \cite{LIMA6, LIMA-JMPA}, in order to prove the sought convergence \eqref{A2}.

\subsection{Convergence of the convective term} \label{ss:conv-conv}

First of all, we notice that, owing to \eqref{w1}, for proving \eqref{A2} it is enough to show that
\begin{equation} \label{A2_bis}
	\int_0^T \intO{ \vc{H}^\perp [\vme] \otimes \vc{H}^\perp [\vme] : \Grad \bfphi } \dt \to 0 \qquad \mbox{ as }\quad \ep \to 0
	\end{equation}
for any test function $\bfphi$ satisfying \eqref{cond:test-f}. Next, omitting an approximation procedure, based
on the spectral decomposition of the wave operator (see the details in \tsl{e.g.} \cite[Sections 5.4.5, 5.4.6]{FeNo6A}),
we may assume that all the quantities appearing in \eqref{A2_bis} and in the wave system \eqref{eq:wave}
are smooth with respect to the space variable.

In light of the previous consideration, we can perform an integration by parts in \eqref{A2_bis} and then compute
\begin{align*}
\Div( \vc{H}^\perp [\vme] \otimes \vc{H}^\perp [\vme] ) &=  \Div(\vc{H}^\perp [\vme] )\ \vc{H}^\perp [\vme] + \vc{H}^\perp [\vme]\cdot\Grad \vc{H}^\perp [\vme] \\
&= \Div\vme \ \vc{H}^\perp [\vme] + \frac{1}{2} \Grad \left| \vc{H}^\perp [\vme] \right|^2 + \Curl(\vc{H}^\perp [\vme]) \times \vc{H}^\perp [\vme]\,.
\end{align*}
Observe that the gradient term identically vanishes, whenever integrated against a test function $\bfphi$ as in \eqref{cond:test-f}.
Similarly, $\Curl(\vc{H}^\perp [\vme]) \equiv 0$, as $\vc{H}^\perp [\vme]$ is a perfect gradient.
Therefore, it remains us to deal with the first term appearing in the right--hand side of the previous relation: for this, we use the wave system
\eqref{eq:wave} to write
\begin{align*}
\Div\vme \ \vc{H}^\perp [\vme] &= -\ep\partial_tZ_\ep\ \vc{H}^\perp [\vme] + \ep \Div\vc{h}^1_\ep\ \vc{H}^\perp [\vme] + \ep^\g h^2_\ep\ \vc{H}^\perp [\vme].
\end{align*}
Owing to the presence of the small factors $\ep$ and $\ep^\g$ and the smoothness of all the involved quantities, it is clear that
the last two terms on the right do not contribute to the limit \eqref{A2_bis}, in the sense that
\[
\int_0^T \intO{ \left( \ep \Div\vc{h}^1_\ep\ \vc{H}^\perp [\vme] + \ep^\g h^2_\ep\ \vc{H}^\perp [\vme] \right) \cdot  \bfphi } \dt \to 0 \qquad \mbox{ as }\quad \ep \to 0.
\]
Finally, we use the second equation in \eqref{eq:wave} to compute
\begin{align*}
-\ep\partial_tZ_\ep\ \vc{H}^\perp [\vme] &= - \veps \partial_t\left(Z_\ep\ \vc{H}^\perp [\vme]\right) + \veps Z_\ep\ \partial_t\vc{H}^\perp [\vme]  \br
&= - \veps \partial_t\left(Z_\ep\ \vc{H}^\perp [\vme]\right) +
\veps^\g Z_\ep\left(\mathbb{H}^3_\ep + \vc{h}^4_\ep\right) - \o Z_\ep \Grad Z_\ep\,.
\end{align*}
In particular, since $Z_\ep \Grad Z_\ep = \frac{1}{2} \Grad(Z_\ep)^2$, the previous computations show that the convergence property \eqref{A2_bis}
holds true, so in turn we have proved \eqref{A2}.

The proof of Theorem \ref{TM1} is now completed.


\end{document}